\documentclass{amsart}

\usepackage{amsfonts,amsthm,latexsym,amsmath,amssymb,amscd,amsmath, mathrsfs, epsf}
\usepackage{graphicx}
\newtheorem{theorem}{Theorem}
\newtheorem{proposition}[theorem]{Proposition}
\newtheorem{lemma}[theorem]{Lemma}

\newtheorem{corollary}{Corollary}

\newtheorem{problem}{Problem}

\newcommand{\HP}{\mathcal{HP}}
\newcommand{\bR}{\mathbb{R}}
\newcommand{\bC}{\mathbb{C}}
\newcommand{\ee}{\end{equation}}
\newcommand {\al}{\alpha}
\newcommand {\la}{\lambda}
\newcommand{\be}{\beta}
\newcommand {\ga}{\gamma}

\newcommand{\A}{\mathcal A}
\newcommand{\B}{\mathcal B}

\begin{document}
\title[Multiplier sequences and logarithmic mesh]{Multiplier sequences and logarithmic mesh}

\author[O.~Katkova]{Olga Katkova}

\address{Department of Mechanics \& Mathematics, Kharkov National University, 4 Svobody Sq., Kharkov, 61077, Ukraine}
\email{olga.m.katkova@univer.kharkov.ua}

\author[B. Shapiro]{Boris Shapiro}

\address{   Department of Mathematics,
            Stockholm University,
            S-10691, Stockholm, Sweden}
\email{shapiro@math.su.se}

\author[A.~Vishnyakova]{Anna Vishnyakova}

\address{Department of Mechanics \& Mathematics, Kharkov National University, 4 Svobody Sq., Kharkov, 61077, Ukraine}
\email{anna.m.vishnyakova@univer.kharkov.ua}

\begin{abstract}
In this note we prove a new  result about (finite) multiplier sequences, i.e. linear operators acting diagonally in the standard  monomial basis of $\bR[x]$ and sending polynomials with all real roots to polynomials with all real roots.  Namely, we show that any such operator does not decrease the logarithmic mesh when acting on an arbitrary polynomial having  all  roots real and of the same sign. The logarithmic mesh of such a polynomial is defined as the minimal quotient of its consecutive roots taken in the non-decreasing order. 
\end{abstract}

\maketitle

\section{Introduction}

Denote by $\HP\subset \mathbb {R}[x]$ the set of all real-rooted  (also referred to as {\em hyperbolic}) polynomials. A linear operator $T:\bR[x]\to\bR[x]$ is called a {\em real rootedness preserver} or a 
{\em hyperbolicity preserver} if it preserves $\HP$. (We will also use the short-hand 'HPO' for such operators. A characterization of hyperbolicity preservers was recently obtained in \cite{BB}.)
Given a real-rooted polynomial $p(x)\in \HP$ denote by $M(p)$ its mesh, i.e. the minimal distance between its real roots. (If $p(x)$ has a double real root then $M(p)=0$.) 

One of (very  few known) results about linear operators not decreasing mesh is presented below and is  due originally to M.~Riesz but was written down by A.~Stoyanoff, see \cite{Sto}. 

\begin{proposition}\label{th:Riesz}  For any hyperbolic polynomial $p$ and any real $\lambda$ one has
$$M(p-\lambda p') \ge M(p), $$
where (as above) $M(q)$ denotes the mesh of the polynomial $q$. 
\end{proposition} 

Recall that the well-known Hermite-Poulain theorem, (see \cite{obr}, p.~4) claims that a finite order linear  differential operator  
$T=a_0+a_1\frac{d}{dx}+...+a_k\frac{d^k}{dx^k}$ with constant coefficients is a hyperbolicity preserver  iff its symbol polynomial $Q(t)=a_0+a_1t+...+a_kt^k$ is hyperbolic.

Thus,  Proposition~\ref{th:Riesz} combined with the Hermite-Poulain theorem immediately imply the following statement.  

\begin{theorem}\label{cor:const} 
Any hyperbolicity preserving differential operator with constant coefficients does not decrease the mesh of hyperbolic polynomials. 
\end{theorem}

For the sake of completeness and due to the fact that \cite{Sto} Êis hardly available nowadays we reprove Theorem~\ref{cor:const} below. Our main interest in this note is to find an  analog of Proposition~\ref{th:Riesz} and Theorem~\ref{cor:const}  for another famous class of HPO, namely, for the so-called multiplier sequences  characterized by G.~P\'olya and I.~Schur in \cite {PolyaSchur}.  The relevant basic  notions are as follows. 

Given a sequence $\A=\{\al_n\}_{n=0}^\infty$ of real or complex numbers we denote by $T_\A$ the linear operator acting diagonally in the monomial basic of $\bC[x]$ by 
$$T_\A(x^i)=\al_ix^i.$$
We refer to $T_\A$ as the {\em diagonal operator corresponding to $\A$}.

 We call a sequence $\A=\{\al_n\}_{n=0}^\infty$ of real numbers  a {\em multiplier sequence of the 1st kind},  if its diagonal operator $T_\A$   preserves $\HP$, i.e. sends an arbitrary real rooted polynomial to a real rooted polynomial. The above sequence $\A$ is called a  {\em multiplier sequence of the 2nd kind},  if the above $T_\A$ sends an arbitrary real rooted polynomial whose roots are all of the same sign to a real rooted polynomial.

The following fundamental criterion  was found  in \cite {PolyaSchur}. 

\begin{theorem}\label{th:PSorig}
Let   $\A=\{\al_n\}_{n=0}^\infty$ be a sequence of real numbers
and let $T_\A : \bR[x] \to  \bR[x] $ be the corresponding diagonal operator. Define  $\Phi(t)$
to be the formal power series
 $$\Phi(t) =\sum_{n=0}^\infty\frac{\al_n}{n!}t^n.$$
The following assertions are equivalent:
\begin{itemize}
\item [(i)]  $\A$ is a multiplier sequence of the first kind,

\item [(ii)]  $\Phi(t)$ defines an entire function which is the limit, uniformly on compact
sets, of polynomials with only real zeros of the same sign,

\item [(iii)] Either  $\Phi(t)$ or  $\Phi(-t)$ is an entire function that can be written as
 $$\Phi(t) = Ct^ke^{at}\Pi_{n=1}^\infty(1 +  \ga_n t);$$
where $k\in \mathbb N$, $C\in\bR,$ $a; \ga_n\ge    0$ for all $n\in \mathbb N$ and
$\sum_{n=1}^\infty \ga_n <\infty$,

\item [(iv)] For all nonnegative integers $k$ the polynomial $T[(1+x)^k]$ is hyperbolic with
all zeros of the same sign.
\end{itemize}
\end{theorem}

  Let us also recall the following analog of the 
      P\'olya-Schur theorem in finite degrees. 
         
Consider a diagonal (in the monomial basis) 
      operator $T: \bR_k[x]\to \bR_k[x]$ acting by multiplication of  $x^j$ by
      $\ga_{j},\;j=0,\ldots,k$. A diagonal hyperbolicity preserver  $T$ will be referred to as a {\em multiplier sequence of length $k+1$}
or simply a {\em finite multiplier sequence}. Denoted by $M_{k}\subset \bR^{k+1}$ the set of all finite
           multiplier sequences of length $k+1$.

The following result was originally proved in \cite{CC2}, Theorem 3.7, see 
also \cite {CC1}, Theorem 3.1.  

\begin{theorem}\label{th:MS}
For $T\in \mathbb{GL}(\bR_k[x])$ the following two conditions are equivalent:
\begin{enumerate}
\item[(i)] $T$ is a multiplier sequence of length $k+1$;
\item[(ii)] The polynomial $Q_{T}(t)=\sum_{j=0}^k{\binom{k}{j}}\ga_{j}t^j$
has all real zeros of the same sign.
\end{enumerate}
\end{theorem}

One can identify the semigroup of all finite multiplier sequences of length $k+1$  with the set of  polynomials of degree $k$ having all real roots of the same sign. The usual multiplication of diagonal matrices then corresponds to the so-called Schur-Szeg\"o multiplication of polynomials, see e.g. \cite{Sz}. Namely, the {\em Schur-Szeg\"o product}  $P*Q$ of two polynomials $P(x)=\sum_{j=0}^k\binom{k}{j}a_jx^j$ and $Q(x)=\sum_{j=0}^k\binom{k}{j}b_jx^j$ equals 
$$P*Q=\sum_{j=0}^k\binom{k}{j}a_jb_jx^j.$$ 

Given a polynomial $P$ of degree $k$  with all real roots of the same sign order these roots as $|x_1|\le |x_2|\le ...\le|x_k|$ and define the {\em logarithmic mesh} of $P$ as $$\rm{lmesh}(P)=\min_{j=1,...,k-1} \frac{|x_{j+1}|}{|x_j|}.$$
Obviously, for any polynomial with all real roots of the same sign one has  that $\rm{lmesh}(P)\ge 1$ and $\rm{lmesh}(P)= 1$ if and only if  $P$ has a multiple (real) root. We can  now formulate the new results of this note. An analog of Proposition~\ref{th:Riesz} is as follows.  

\begin{proposition}\label{th:OA} For any $\la>0$  the differential operator $T(p(x))=\la p+xp'$ has the property 
$$\rm{lmesh}(T(p))\ge \rm{lmesh}(p),$$
where $p$ is an arbitrary polynomial  with all real roots of the same sign.
\end{proposition}Ê

The latter proposition can be generalized to the Schur-Szeg\"o product  of two polynomials of the same degree. Namely,    the following statement holds.  

\begin{theorem}\label{th:SchSz} Given two polynomials $P$ and $Q$ of degree $k$  with all  roots of the same sign one has 
$$\rm{lmesh}(P*Q)\ge \max (\rm{lmesh}(P),\rm{lmesh}(Q)).$$
\end{theorem} 

\noindent
{\bf Remark.}
Notice that Theorem~\ref{th:SchSz} Êdoes not follow from Proposition~\ref{th:OA} (contrary to the case of constant coefficients) since not all (finite) multiplier sequences considered as differential operators can be represented as the product of   operators of the form $\la p+xp'$. On the other hand, Proposition~~\ref{th:OA} follows from Theorem~\ref{th:SchSz}  since for any given positive integer $k$ one can represent 
the action of the operator $\la + x\frac{d}{dx}$ on polynomials of degree $k$ as the Schur-Szeg\"o composition with  an appropriate polynomial, see \S~2. Also we want to mention that the original proof of Proposition~\ref{th:Riesz} by M.~Riesz can not be extended to multiplier sequences and logarithmic mesh. This part of the present paper is new. 



\medskip
\noindent
{\it Acknowledgements.} The authors are grateful to P.~Br\"anden and V.~Kostov for numerous discussions of the area related to HPO. The first and the third  authors want to acknowledge the hospitality of Department of Mathematics, Stockholm  University during their visit to Stockholm in  March-April 2010 supported by a personal grant of Professor Br\"anden.  

\medskip

\section{Proofs}Ê

In the proofs of  Theorems~\ref{cor:const}  and \ref{th:SchSz} we will use the following three facts  collected in a lemma below.

\begin{lemma}\label{lm:all}

\begin{itemize}

\item[(i)]  Given two real polynomials $f$ and $g$ of the same degree one has that the pencil  $cf(x)+dg(x)$ consists of hyperbolic polynomials if and only if 
  $f$ and  $g$ has all  real and (non-strictly) interlacing roots. (This is known under the name Obreschkov's theorem, see \cite{obr} although it has been (re)discovered many times by different authors in the past.)

\item[(ii)] Given two polynomials $P$ and $Q$ of the same degree satisfying the conditions that $P$ and $Q$ are hyperbolic and, additionally, all roots of $Q$ are of the same sign then 
their Schur-Szeg\"o composition $(P*Q)$ is hyperbolic. (This is a special case of the well-known Malo-Schur-Szeg\"o theorem, see e.g. \cite{Sz}).

\item[(iii)] Take $P(x)=a(x+x_1)(x+x_2)\cdot \ldots \cdot (x+x_k),$
where $0< x_1 < x_2 < \ldots < x_k ,$ and choose $\lambda > 1$. Then
the zeros of $P(x)$ and $P(\lambda x)$ are interlacing if and only
if $\lambda < \rm{lmesh}(P)$.

\end{itemize}

\end{lemma}

Let us settle Theorem~~\ref{cor:const}. 

\begin{proof}  Given an linear ordinary  differential operator $A$ with constant coefficients and of finite order suppose that there exists $P \in \mathbb{R}[x]$ with
all real zeros such that $ \rm{mesh}( A(P)) < \rm{mesh}(P)$. Choose $\lambda$ satisfying the inequalities $0 \leq \rm{mesh}(A(P)) < \lambda <
\rm{mesh}(P)$. Then by (iii) the zeros of polynomials $A(P)(x)$ and
$A(P)(x+\lambda )$ are not interlacing. Using  (i) we get that  there exist $c,
d \in \mathbb{R}$ such that the polynomial
$cA(P)(x)+dA(P)(x+\lambda )$ has a nonreal zero. But
$cA(P)(x)+dA(P)(x+\lambda ) = A(L)(x)$, where $L(x)=
cP(x)+dP(x+\lambda )$. Since $\lambda < \rm{mesh}(P)$ the zeros of
$P(x)$ and $P(x+\lambda )$ are interlacing, and by (i) we conclude 
that all zeros of $L$ are real. But, this implies that  all zeros of $A(L)$ should be real as well. This contradiction finishes the proof. 
\end{proof}

We now settle Theorem~\ref{th:SchSz}. 

\begin{proof}ÊGiven $P$ and $Q$ two polynomials of the same degree with all real roots of the same sign consider $S(x):= (P*Q)(x)$. Assume that
$\rm{lmesh}(S) < \rm{lmesh}(P)$ and choose $\lambda$ such that  $1\leq
\rm{lmesh}(S) < \lambda < \rm{lmesh}(P)$. By (iii) since
$\rm{lmesh}(S) < \lambda$ the zeros of polynomials $S(x)$ and
$S(\lambda x)$ are not interlacing. By (i) there exist $c, d
\in \mathbb{R}$ such that the polynomial $cS(x)+dS(\lambda x)$ has
a nonreal zero. We have
$$cS(x)+dS(\lambda x)=\sum_{j=0}^k {k \choose j}  (ca_j +d\lambda^j a_j) b_j x^j .$$
Denote by $L(x)= \sum_{j=0}^j {k \choose j}  (ca_j +d\lambda^j a_j)  x^j
=cP(x) +d P(\lambda x)$. Then one gets 
$$cS(x)+dS(\lambda x)= (L*Q)(x).$$
By (iii) the inequality  $\lambda < \rm{lmesh}(P)$ implies that the zeros of polynomials
$P(x)$ and $P(\lambda x)$ are interlacing. Then by (i) all the 
zeros of $L(x)$ are real, and all the  zeros of $Q(x)$ are real and of
the same sign. Then by (ii) all the zeros of $cS(x)+dS(\lambda x)$ should be 
real. This contradiction finishes the proof. 
\end{proof} 

Let us finally deduce Proposition~\ref{th:OA} from Theorem~\ref{th:SchSz}. 

\begin{proof} One can easily check that for any two polynomials $P$ and $Q$ of the same degree the relation $(P+axP')* Q=(P*Q)+ax(P*Q)'$ holds. When $P(x)=(x+1)^k$ one gets
$$Q+axQ'=((x+1)^{k-1}((1+ak)x+1))*Q.$$
Therefore, for $a>0$ the action of the differential operator $1+ax\frac{d}{dx}$ on the polynomial $Q$ coincides with its composition with the polynomial with all negative roots and the result follows. 
\end{proof}

\section{Final remarks} 

\noindent
{\bf 1.}Ê Theorem~~\ref{th:SchSz} shows that for any two finite multiplier sequence of the same length the logarithmic mesh of their composition is greater than or equal to the maximum of their logarithmic meshes.  One can try to generalize this result to the case of usual (infinite) multiplier sequences. Given a multiplier sequence $\A=\{\al_n\}_{n=0}^\infty$ define its logaritmic mesh as 
$$\rm{lmesh}(\A)=\inf_{k=0,...,\infty}\rm{lmesh}(\A_k),$$ 
where $\A_k$ is the $k$-th truncation of $\A$, i.e. its initial finite segment of length $k+1$. Then Theorem~~\ref{th:SchSz} immediately implies the following.  
\begin{corollary}\label{cor:inf} For any infinite multiplier sequences $\A$ and $\B$ one has 
$$\rm{lmesh}(\A\B)\ge \max(\rm{lmesh}(\A),\rm{lmesh}(\B)).$$
\end{corollary}

\begin{problem} Describe the class of multiplier sequences whose logarithmic mesh is strictly greater than $1$. 
\end{problem} 

\noindent
{\bf 2.} ÊTheorems~~\ref{cor:const} and \ref{th:SchSz}  are examples of a statement that an appropriate version of mesh on the space of hyperbolic polynomials is preserved under the action of the corresponding class of hyperbolicity preservers. At the moment we have no idea what kind of mesh one should associate to a arbitrary HPO so that its action on hyperbolic polynomials will increase it for generic hyperbolic polynomials. Such a notion would be highly desirable.

%


\begin{thebibliography}{8}

\bibitem{BB}  J.~Borcea, P.~Br\"and\'en,   {\em P\'olya-Schur master theorems for circular domains and their boundaries.} Ann. of Math. (2) 170 (2009), no. 1, 465--492.

    

\bibitem{CC1} 
T.~Craven, G.~Csordas, {\em Problems and theorems in the theory of multiplier sequences}, Serdica Math. J. {\bf 22} (1996), 515--524.


\bibitem {CC2} T.~Craven and G.~Csordas, 
{\em Multiplier sequences for fields.} Illinois J. Math. 21 (1977), no. 4, 801--817.






\bibitem {La}ÊE.~Laguerre, Sur quelques points de la th\'eorie des  \'equations num\'eriques, Acta Math. 4 (1884), 97--120. 

    
    \bibitem{obr}
   N.~Obreschkov, Verteilung und Berechnung der Nullstellen reeller Polynome, VEB Deutscher
Verlag der Wissenschaften, Berlin, 1963.

\bibitem {PolyaSchur} G.~P\'olya and J.~Schur, {\em \"Uber zwei Arten von Faktorenfolgen in der Theorie der algebraische Gleichungen}, J. Reine und Angew. Math., v. 144 (1914), 89-113. 

\bibitem{Sto} A. Stoyanoff,  {\em Sur un theoreme de M Marcel Riesz,}  Nouvelles Annales de Mathematique 1
(1926), 97--99. 


\bibitem{Sz} G.~Szeg\"o, Bemerkungen zu einem Satz von 
J.~H.~Grace \"uber die Wurzeln algebraischer Gleichungen, 
Math.Z., (2) vol. 13, (1922), 28-55.  

\end{thebibliography}
\end{document}